\documentclass[11pt,reqno]{amsart}

\usepackage[text={160mm,240mm},centering]{geometry}            % See geometry.pdf to learn the layout options. There are lots.
\usepackage{amssymb,amsmath,setspace,geometry,indentfirst,changepage,inputenc,amsthm}
\usepackage{mathrsfs,amsfonts,savesym,graphicx,bm,color}

\usepackage{graphicx}
\usepackage{amssymb}
\usepackage{epstopdf}
\usepackage{dsfont}
\usepackage{placeins}
\usepackage{pifont}
\usepackage{romannum}

\usepackage{lscape}

\usepackage{tikz}
\usetikzlibrary{positioning}
\usepackage{xcolor}

\usepackage{float}

\usepackage[all,2cell]{xy}    % ÒýÈë xy °ü
\UseAllTwocells               % ²¢ÇÒÐèÒª»­ 2-cell

\usepackage{amsmath,amsthm}

\usepackage[mathscr]{eucal}
\usepackage[all]{xy}
\usepackage{mathrsfs}
\usepackage{hyperref}
\usepackage{color}

\DeclareMathOperator{\Spec}{Spec}
\DeclareMathOperator{\Hom}{Hom}
\DeclareMathOperator{\End}{End}
\DeclareMathOperator{\Ima}{Im}

\newtheorem{theorem}{Theorem}[section]
\newtheorem{lemma}[theorem]{Lemma}

\newtheorem*{main thm}{Main Theorem}

\newtheorem{conj}[theorem]{Conjecture}
\newtheorem{cor}[theorem]{Corollary}
\newtheorem{prop}[theorem]{Proposition}

\newtheorem{def-lemma}[theorem]{Definition-Lemma}
\newtheorem{def-thm}[theorem]{Definition-Theorem}

\newtheorem{thm}[theorem]{Theorem}

\theoremstyle{definition}

\newtheorem{definition}[theorem]{Definition}

\newtheorem{remark}[theorem]{Remark}

\allowdisplaybreaks

\title[]{Hilbert series of jet schemes of determinantal varieties}

\author{Yongxin Xu}
\address{
	Department of Mathematical Sciences,
	Tsinghua University,
	Beijing, 100084, P. R. China.}
\email{xuyongxi22@mails.tsinghua.edu.cn}

\author{Huaiqing Zuo}
\address{Department of Mathematical Sciences,
	Tsinghua University,
	Beijing, 100084, P. R. China.}
\email{hqzuo@mail.tsinghua.edu.cn}

\begin{document}

	\setcounter{page}{1}
	\pagenumbering{arabic}
	
	\maketitle
	
	\parskip=.3em
	
	\setcounter{tocdepth}{1}
	%\tableofcontents

	\begin{abstract}
		Let $X$ be the variety of $m\times n$ matrices of rank less than $r$ over an algebraically closed field $K$ of characteristic zero. We calculate the Hilbert series of the principal component of the $k$th jet scheme of $X$. This solves the Hilbert series assertion of a conjecture of Chen and Zuo. The proof uses a one-parameter deformation of $K[t]/(t^{k+1})$ and classical invariant theory for the general linear group.
		
		
	\end{abstract}
	
	%\begin{document}
	%	\maketitle
	\section{Introduction}\label{s1}
	
	\par The theory of jet schemes and arc spaces was introduced into singularity theory by Nash, with the aim of understanding resolutions of singularities through families of arcs passing through singular points, see \cite{c1}. Ein and Mustaţă studied some invariants of singularities using jet schemes, see \cite{c2}. Determinantal varieties form an important class of examples in singularity theory. This paper studies the Hilbert series of jet schemes of determinantal varieties.
	
	\par Throughout this paper, $K$ denotes an algebraically closed field of characteristic zero. Let $r, m, n\in \mathbb{Z}_{>0}$ be integers satisfying $r\le m\le n$. Denote by $M(m\times n, K)$ the space of $m\times n$ matrices with entries in $K$. It is naturally isomorphic to the affine space $\mathbb{A}_{K}^{mn}$. For a fixed $m \times n$ matrix $(x_{i,j})_{1\le i\le m, 1\le j\le n}$, all $r \times r$ minors of the matrix generate an ideal $I^{m,n}_{r}$ in $K[x_{i,j}|1\le i\le m, 1\le j\le n]$, which defines an algebraic variety $\mathcal{L}^{m,n}_{r}$. We say $\mathcal{L}^{m,n}_{r}$ is a determinantal variety of $m \times n$ matrix.
	
	\par Let $k\in \mathbb{Z}_{\ge 0}$. We write $J_{k}(\mathcal{L}_{r}^{m,n})$ for the $k$th jet scheme of $\mathcal{L}^{m,n}_{r}$. The precise definition is given in Section 2.1. As is immediate from that definition, $J_{0}(\mathcal{L}_{r}^{m,n})$ is isomorphic to $\mathcal{L}_{r}^{m,n}$, and $J_{k}(\mathcal{L}_{r}^{m,n})$ is defined by a homogeneous ideal. Hence we can discuss their Hilbert series.
	
	\par Conca and Herzog in \cite{c3} calculate the Hilbert series of classical determinantal varieties.
	
	\begin{thm}\label{t1.1}\textnormal{(\cite{c3}, Corollary 1)}
		The Hilbert series of $\mathcal{L}^{m,n}_{r+1}$ is
		\begin{equation}\label{eq1}
			\centering
			H_{\mathcal{L}^{m,n}_{r+1}}(z)= \frac{\det\big( \underset{k\ge 0}{\sum}\binom{m-i}{k} \binom{n-j}{k} z^{k}\big)_{1\le i\le r,1\le j\le r}}{z^{\binom{r}{2} }(1-z)^{r(m+n-r)}}.
		\end{equation}
	\end{thm}
	
	\par It is complicated to compute the Hilbert series for general $J_{k}(\mathcal{L}_{r}^{m,n})$. The calculation of Hilbert series of $J_{1}(\mathcal{L}_{2}^{m,n})$ was finished by three papers \cite{c4}, \cite{c5}, and \cite{c6}, using Gröbner basis and combinatorial methods. The result is as follows:
	\begin{equation}\label{eq2}
		\centering
		H_{J_{1}(\mathcal{L}^{m,n}_{2})}(z)=\left\{
		\begin{array}{rcl}
			& H_{\mathcal{L}^{m,n}_{2}}(z)^{2}, & \quad m=2, \\
			& H_{\mathcal{L}^{m,n}_{2}}(z)^{2}+\frac{1}{(1-z)^{mn}}-H_{\mathcal{L}^{m,n}_{3}}, & \quad m\ge 3.
		\end{array}
		\right.
	\end{equation}
	Using a similar approach, Chen and Zuo in \cite{c7} computed the Hilbert series of $J_{1}(\mathcal{L}_{3}^{3,n})$, and Chen, Xu, and Zuo in \cite{c8} computed the Hilbert series of $J_{2}(\mathcal{L}_{2}^{2,n})$. Their results are as follows:
	
	\begin{theorem}\label{t1.2}\textnormal{(\cite{c7}, Theorem B)}
		The Hilbert series of $J_{1}(\mathcal{L}_{3}^{3,n})$ is
		\begin{equation}\label{eq3}
			\centering
			H_{J_{1}(\mathcal{L}_{3}^{3,n})}(z)=\big(  \frac{1+(n-2)z+\frac{(n-1)(n-2)}{2}z^{2}}{(1-z)^{2n+2}}\big)^{2}.
		\end{equation}
	\end{theorem}
	
	\begin{theorem}\label{t1.3}\textnormal{(\cite{c8}, Theorem 6.1)}
		The Hilbert series of $J_{2}(\mathcal{L}_{2}^{2,n})$ is
		\begin{equation}\label{eq4}
			\centering
			H_{J_{2}(\mathcal{L}^{2,n}_{2})}(z)=\big( \frac{1+(n-1)z}{(1-z)^{n+1}}\big)^{3}.
		\end{equation}
	\end{theorem}
	
	\par Combining (\ref{eq1})-(\ref{eq4}), it is natural to ask whether there is a relation between $H_{J_{k}(\mathcal{L}_{r}^{m,n})}(z)$ and $H_{\mathcal{L}^{m,n}_{r}}(z)$. This question should be concerned with principal component of $J_{k}(\mathcal{L}_{r}^{m,n})$, rather than $J_{k}(\mathcal{L}_{r}^{m,n})$ itself. We denote by $P_{k}(\mathcal{L}_{r}^{m,n})$ the principal component of $J_{k}(\mathcal{L}_{r}^{m,n})$. The precise definition is given in Section 2.1. In \cite{c6}, the authors obtained the following result:
	\begin{equation*}
		H_{P_{1}(\mathcal{L}^{m,n}_{2})}(z)=H_{\mathcal{L}^{m,n}_{2}}(z)^{2}.
	\end{equation*}
	All jet schemes of $\mathcal{L}_{m}^{m,n}$ are irreducible, see Theorem 3.1 in \cite{c15}. Hence we have $J_{1}(\mathcal{L}_{3}^{3,n})=P_{1}(\mathcal{L}_{3}^{3,n})$ and $J_{2}(\mathcal{L}_{2}^{2,n})=P_{2}(\mathcal{L}_{2}^{2,n})$. Motivated by these results, Chen and Zuo in \cite{c7} proposed following conjecture.
	
	\begin{conj}\label{main-conj}\textnormal{(\cite{c7}, Conjecture 1.2)}
		The Hilbert series of the principal component of $k$th jet scheme of $\mathcal{L}_{r}^{m,n}$ is exactly the $(k+1)$th power of the Hilbert series of $\mathcal{L}_{r}^{m,n}$.
	\end{conj}
	
	\par As said in \cite{c7}, there are a lot of examples showing that the conjecture still
	holds for the general case. We prove the conjecture in this paper.
	
	\begin{main thm}\label{mt}
		The Hilbert series of $P_{k}(\mathcal{L}_{r}^{m,n})$ is exactly the $(k+1)$th square of the Hilbert series of $\mathcal{L}_{r}^{m,n}$.
	\end{main thm}
	
	\begin{remark}\label{r1.1}
		A technical difficulty in earlier calculations of $H_{P_{k}(\mathcal{L}_{r}^{m,n})}(z)$ is that Gröbner basis calculations become increasingly complicated as the order of the jet scheme grows. Our approach is different and uses tools from representation theory.
	\end{remark}
	
	\par This paper is organized as follows.
	
	\par In Section \ref{s2}, we provide preliminaries about jet schemes, Hilbert series, and representation theory. In Section \ref{s3}, we prove the main theorem.

	\section{Preliminaries}\label{s2}
	
	\subsection{Jet scheme and principal component.}
	
	\par $\quad$
	
	\par Let $X$ be a scheme of finite type over $K$. For all $k\in \mathbb{Z}_{\ge 0}$, we say a $K$-scheme $J_{k}(X)$ is the $k$th jet scheme of $X$, if for every $K$-algebra $A$, there exists a  functorial bijection
	\begin{align*}
		\Hom (\Spec A,J_{k}(X))\simeq \Hom (\Spec A[t]/(t^{k+1}),X).
	\end{align*}
	The $k$th jet scheme always exists, and it is also a scheme of finite type over $K$. It is unique up to a canonical isomorphism. One can see that $J_{0}(X)$ is isomorphic to $X$. We define the truncation map $\pi_{k}:J_{k}(X)\to X$ via $A[t]/(t^{k+1})\to A[t]/(t)$. Let $f:X\to Y$ be a morphism between $K$-schemes. Then it naturally induces $J_{k}(f):J_{k}(X)\to J_{k}(Y)$ via the bijection $\Hom (\Spec A,J_{k}(X))\simeq \Hom (\Spec A[t]/(t^{k+1}),X)$. We refer the reader to \cite{c2} for more details.
	
	\par For an affine variety, one can explicitly describe its jet scheme as well as the morphisms between jet schemes. Suppose $X=\Spec K[x_{1},\cdots ,x_{m}]/I$ is an affine variety in $\mathbb{A}_{K}^{m}$ defined by the ideal $I$. We write $x_{i}(t)=\sum\limits_{j=0}^{\infty}x_{i}^{(j)}t^{j}$, where $t$ is an indeterminate. The $k$th jet scheme of $X$ is $J_{k}(X)=\Spec K[x_{i}^{(j)}|1\le i\le m,0\le j\le k]/J$, where $J$ is generated by
	\begin{align*}
		\{f^{(0)},\cdots ,f^{(k)}|f(x_{1}(t),\cdots ,x_{m}(t))\equiv \sum\limits_{l=0}^{k}f^{(l)}(x_{i}^{(j)})t^{l}\text{ (mod $t^{k+1}$) for some $f\in I$}\}.
	\end{align*}
	By taking $x_{i}^{(j)}=0$ for all $j\ge 1$, we can define the truncation map $\pi_{k}:J_{k}(X)\to X$. In these coordinates, $\pi_k$ sends a jet to its constant term.
	
	\par Let $X=\Spec K[x_{1},\cdots ,x_{m}]/I_{1}$ and $Y=\Spec K[y_{1},\cdots ,y_{n}]/I_{2}$ be two affine varieties. Suppose $f=(f_{1},\cdots ,f_{n}):X\to Y$ is a map with $f_{1},\cdots ,f_{n}$ regular functions. We can define 
	\begin{align*}
		J_{k}(f)=(f_{1}^{(0)},\cdots ,f_{1}^{(k)},\cdots ,f_{n}^{(0)},\cdots ,f_{n}^{(k)}):J_{k}(X)\to J_{k}(Y),
	\end{align*} 
	where $f_{i}^{(j)}$ is defined via $f_{i}(x_{1}(t),\cdots ,x_{m}(t))\equiv \sum\limits_{l=0}^{k}f_{i}^{(l)}(x_{i}^{(j)})t^{l}$ (mod $t^{k+1}$).
	
	\begin{lemma}\label{l2.1}\textnormal{(\cite{c2}, Lemma 2.3)}
		Let $X$ be a scheme of finite type over $K$, let $k\in \mathbb{Z}_{\ge 0}$, and let $U\subset X$ be an open subset. Then the $k$th jet scheme of $U$ is $J_{k}(U)=\pi_{k}^{-1}(U)$.
	\end{lemma}
	
	\begin{definition}\label{d2.1}
		Let $X$ be an algebraic variety over $K$, let $k\in \mathbb{Z}_{\ge 0}$, and let $\pi_{k}:J_{k}(X)\to X$ be the truncation map of jet scheme. Suppose $U$ is the smooth locus of $X$. The principal component of $J_{k}(X)$ is defined to be $P_{k}(X)=\overline{\pi_{k}^{-1}(U)}\subset J_{k}(X)$.
	\end{definition}	
	
	\par For an affine $K$-variety $X\subset \mathbb{A}_{K}^{m}$, we denote by $p_{k}(X)$ the ideal defining $P_{k}(X)$ in $\mathbb{A}_{K}^{m(k+1)}$. If $X$ is defined by a homogeneous ideal, then $p_{k}(X)$ is homogeneous since $J_{k}(X)$ is defined by a homogeneous ideal. Hence it makes sense to consider the Hilbert series of $P_{k}(\mathcal{L}_{r}^{m,n})$.
	
	\begin{lemma}\label{l2.2}\textnormal{(\cite{c10}, Lemma 2.7)}
		Suppose $X$ is an irreducible $K$-variety, and let $k\in \mathbb{Z}_{\ge 0}$. The principal component $P_{k}(X)$ is an irreducible component of $J_{k}(X)$.
	\end{lemma}
	
	\begin{cor}\label{c2.1}
		Let $X$ and $Y$ be two affine varieties over $K$, let $f:X\to Y$ be a dominant map, and let $k\in \mathbb{Z}_{\ge 0}$. If $X$ is nonsingular and irreducible, then
		\begin{enumerate}
			\item $P_{k}(Y)=\overline{J_{k}(f)\big( J_{k}(X)\big)}$;  
			\item The image of the homomorphism of coordinate rings $J_{k}(f)^{*}:K[J_{k}(Y)]\to K[J_{k}(X)]$ is isomorphic to $K[P_{k}(Y)]$. 
		\end{enumerate}
	\end{cor}
	\begin{proof}[Proof. ]
		Let $W$ be the smooth locus of $Y$. By Theorem 21.6.4(b) in \cite{c11}, there exists a nonempty open subset $U\subset f^{-1}(W)$ such that $f|_{U}$ is smooth. Since smooth maps are open, $V=f(U)$ is a nonempty open subset of $W$, and $f|_{U}:U\to V$ is smooth and surjective. By Lemma 2.4 in \cite{c10}, the induced map $J_{k}(U)\to J_{k}(V)$ is smooth and surjective. 
		\par By Proposition 2.6 in \cite{c10}, $J_{k}(X)$ and $J_{k}(W)$ are irreducible. By Lemma \ref{l2.1}, $J_{k}(U)$ and $J_{k}(V)$ are nonempty open subsets of $J_{k}(X)$ and $J_{k}(W)$, respectively, and hence are dense. Taking closures in $J_{k}(Y)$, we obtain
		\begin{align*}
			\overline{J_{k}(f)\big( J_{k}(X)\big)}=\overline{J_{k}(f)\big( J_{k}(U)\big)}=\overline{J_{k}(V)}=\overline{J_{k}(W)}=P_{k}(Y).
		\end{align*}
		This proves (1). By (1), the kernel of $J_{k}(f)^{*}$ is precisely the ideal defining $P_{k}(Y)$ in $K[J_{k}(Y)]$. The second assertion follows from the first isomorphism theorem.
	\end{proof}

	\subsection{Hilbert series.}
	
	\par $\quad$
	
	\par For a finitely generated graded $K$-algebra $R$, we denote by $H_{R,K}(z)$ the Hilbert series of $R$. If $L/K$ is a field extension, then $H_{R\otimes_{K}L,L}(z)=H_{R,K}(z)$. If no ambiguity arises, we will write $H_{R}(z)$ for short.
	
	\begin{definition}\label{d2.2}
		Let $R=\bigoplus\limits_{d\ge 0}R_{d}$ be a graded $K$-algebra. We say $R$ is standard if $R_{0}=K$, $R=K[R_{1}]$, and $\dim_{K}R_{1}<\infty$.
	\end{definition}
	
	\par One can see that a standard graded $K$-algebra is finitely generated, and the coordinate ring of an affine $K$-variety is a standard graded $K$-algebra. 
	
	\begin{prop}\label{p2.1}\textnormal{(\cite{c9}, Theorem 0.1)}
		Let $R_1$ and $R_2$ be standard graded $K$-algebras. If they are isomorphic as $K$-algebras, then they are isomorphic as graded $K$-algebras.
	\end{prop}
	
	\begin{lemma}\label{l2.3}
		Let $R$ be a graded $K[\varepsilon]$-algebra, where $\varepsilon$ is an indeterminate. Suppose $R_{d}$ is a finitely generated free $K[\varepsilon]$-module for all $d\ge 0$. Then for any $a\in K$, we have $H_{R/(\varepsilon -a)R,K}(z)=H_{R\otimes_{K[\varepsilon]}K(\varepsilon),K(\varepsilon)}(z)$.
	\end{lemma}
	\begin{proof}[Proof. ]
		Set $r_{d}=\operatorname{rank}_{K[\varepsilon]}R_{d}$. Since $R_{d}$ is a finitely generated free $K[\varepsilon]$-module, for every $a\in K$ and $d\ge 0$, we have
		\begin{align*}
			\dim_{K}\big( R_{d}/(\varepsilon-a)R_{d}\big)=r_{d}=\dim_{K(\varepsilon)}\big( R_{d}\otimes_{K[\varepsilon]}K(\varepsilon)\big).
		\end{align*}
		Since $\varepsilon$ has degree zero, these are the dimensions of the degree $d$ components of $R/(\varepsilon-a)R$ and $R\otimes_{K[\varepsilon]}K(\varepsilon)$, respectively. Comparing the coefficients of their Hilbert series gives the result.
	\end{proof}

	\subsection{Flag variety.}
	
	\par $\quad$
	
	\begin{definition}\label{d2.3}
	    Let $V$ be a $K$-vector space of dimension $N$, and let $l\in \mathbb{Z}_{>0}$. Suppose $b_{1},\cdots ,b_{l}\in \mathbb{Z}_{>0}$ be $l$ positive integers satisfying $b_{1}+\cdots +b_{l}=N$. The flag variety is
	    \begin{align*}
	    	F_{b_{1},\cdots ,b_{l}}(V)=\{ (F_{0},\cdots ,F_{l})|0=F_{l}\subset F_{l-1}\subset \cdots \subset F_{0}=V,\dim_{K}F_{i-1}-\dim_{K}F_{i}=b_{i}\}.
	    \end{align*}
	\end{definition}
	
	\par Let $V$ be a $K$-vector space, and let $d\in \mathbb{Z}_{>0}$. We denote by $Gr(d,V)$ the Grassmannian variety defined by $d$-dimensional subspaces of $V$. Following the notation in Definition \ref{d2.3}, the map 
	\begin{align*}
		F_{b_{1},\cdots ,b_{l}}(V)\to \prod\limits_{j=1}^{l-1}Gr(N-b_{1}-\cdots -b_{j},V),(F_{0},\cdots ,F_{l})\mapsto (F_{1},\cdots ,F_{l-1})
	\end{align*}
	induces the $K$-variety structure on flag variety. Flag variety is smooth and irreducible. 
	
	\par We now describe the cotangent space of the flag variety $F=F_{b_{1},\cdots ,b_{l}}(V)$ following the notation in Definition \ref{d2.3}. Let $F_{\bullet}=(F_{0},\cdots ,F_{l})$ be a flag in $F$. We define
	\begin{align*}
		\mathfrak{p}_{F_{\bullet}}=\{ A\in \End(V)|AF_{i}\subset F_{i}\text{ for all $i$}\}\\
		\mathfrak{u}_{F_{\bullet}}=\{ D\in \End(V)|DF_{i-1}\subset F_{i}\text{ for all $i$}\}
	\end{align*}
	The tangent space $T_{F_{\bullet}}F$ is isomorphic to $\End(V)/\mathfrak{p}_{F_{\bullet}}$. We have an isomorphism
	\begin{align*}
		\mathfrak{u}_{F_{\bullet}}\simeq T^{*}_{F_{\bullet}}F,D\mapsto ([A]\mapsto tr(DA)).
	\end{align*}
	The cotangent bundle $T^{*}F$ can be embedded into $\End(V)\times F$, and the isomorphism above shows
	that $T^{*}F=\{ (D,F_{\bullet})\in \End(V)\times F|D\in \mathfrak{u}_{F_{\bullet}}\}$. The projection to $\End{V}$ induces a projective map $\mu_{0}:T^{*}F\to \End(V)$. The left multiplication action of $GL(V)$ on $F$ induces the conjugate action on the cotangent bundle. We denote by 
	\begin{align*}
		\mathcal{N}_{b_{1},\cdots ,b_{l}}=\overline{\{ gJ_{0}g^{-1}|g\in GL(V)\}}\subset \End(V)
	\end{align*}
	the nilpotent orbit of $F$, where $J_{0}$ is the Jordan matrix consisting of nilpotent Jordan blocks of sizes $b_{1},\cdots ,b_{l}$ placed successively along the diagonal. Then we obtain a projective birational surjective map $\mu_{0}:T^{*}F\to \End(V)$, which is the Springer resolution. The variety $\mathcal{N}_{b_{1},\cdots ,b_{l}}$ is normal. For more details, we refer the reader to \cite{c13}.
	
	\par Let $\lambda_{1},\cdots ,\lambda_{l}$ be different elements in $K$. We introduce a certain deformation of the Springer resolution following. We set
	\begin{align*}
		Z=\{ (D,F_{\bullet},\varepsilon)\in \End(V)\times F\times \mathbb{A}_{K}^{1}| (D-\lambda_{i}\varepsilon I)F_{i-1}\subset F_{i} \text{ for all $1\le i\le l$}\}.
	\end{align*}
	This gives a Zariski locally trivial fiber bundle on $F\times \mathbb{A}_{K}^{1}$. The projection from $Z$ to $\mathbb{A}_{K}^{1}$ and $\End(V)\times \mathbb{A}_{K}^{1}$ are denoted by $\pi:Z\to \mathbb{A}_{K}^{1}$ and $\mu:Z\to \End(V)\times \mathbb{A}_{K}^{1}$, respectively. The zero fiber $\pi^{-1}(0)$ is isomorphic to $T^{*}F$. For all $c\in K$ nonzero, the fiber $\pi^{-1}(c)$ is isomorphic to $\{ gJ_{c}g^{-1}|g\in GL(V)\}$, where $J_{c}$ is the block diagonal matrix with diagonal blocks $\lambda_{1}cI_{b_{1}},\cdots ,\lambda_{l}cI_{b_{l}}$ arranged successively along the diagonal.
	
	\begin{prop}\label{p2.4}\textnormal{(\cite{c14}, Lemma 5.1)}
		Notations are as above. The variety $Z$ is smooth and irreducible, and the map $\mu$ is projective.
	\end{prop}

	\subsection{Rational representation of general linear group.}
	
	\par $\quad$
	
	\par In this subsection, we assume that $V$ is a $K$-vector space of finite dimension, and $G=GL(V)$ is an algebraic group.
	
	\begin{definition}\label{d2.4}
		Let $M$ be a $KG$-module. We say $M$ is a rational $G$-module if $M$ is the union of all its finite dimensional submodules and the action of $G$ on each finite dimensional submodule is algebraic.
	\end{definition}
	
	\par Let $X$ be an affine $K$-variety with an algebraic $G$-action. The map $G\times K[X]\to K[X],g\cdot f\mapsto f\circ g^{-1}$ induces a rational $G$-module structure on $K[X]$. For a rational $G$-module $M$, we write $M^{G}=\{ v\in M|gv=v\text{ for all $g\in G$}\}$.
	
	\begin{prop}\label{p2.2}
		Suppose 
		\begin{align*}
			0\to M_{1}\to M_{2}\to M_{3}\to 0
		\end{align*}
		is a short exact sequence of rational $G$-modules and $G$-equivariant $K$-linear maps. Then 
		\begin{align*}
			0\to M_{1}^{G}\to M_{2}^{G}\to M_{3}^{G}\to 0
		\end{align*}
		is exact.
	\end{prop}
	
	\begin{prop}\label{p2.3}
		Let $r_{1},r_{2}>0$ be two positive integers, and let $W=(V^{*})^{r_{1}}\oplus V^{r_{2}}$. Suppose $\End(V)\oplus W$ has the $G$-action
		\begin{align*}
			g\cdot (X,u_{1},\cdots ,u_{r_{1}},v_{1},\cdots ,v_{r_{2}})=(gXg^{-1},u_{1}\circ g^{-1},\cdots ,u_{r_{1}}\circ g^{-1},gv_{1},\cdots ,gv_{r_{2}}).
		\end{align*}
		Then we have $K[\End(V)\oplus W]^{G}=K[tr(X^{q}),u_{i}(X^{l}v_{j})|q\ge 1,l\ge 0,1\le i\le r_{1},1\le j\le r_{2}]$.
	\end{prop}
	\begin{proof}[Proof. ]
		Fix a basis of $V$, and regard $u_{i}$ as a row vector and $v_{j}$ as a column vector. Put $Y_{ij}=v_{j}u_{i}$. For $z\in K^{*}$, the scalar matrix $zI$ fixes $X$ and sends $(u_{i},v_{j})$ to $(z^{-1}u_{i},zv_{j})$. Thus in every monomial of an invariant polynomial the number of vector coordinates equals the number of covector coordinates, counted with multiplicity. Indeed, a monomial with $p$ vector coordinates and $q$ covector coordinates is multiplied by $z^{p-q}$, and its coefficient must agree for every $z\in K^{*}$. Pairing these coordinates shows that every invariant is a polynomial in the entries of $X$ and the $Y_{ij}$, since each paired product is an entry of some $v_{j}u_{i}$.
		\par Let $B$ be the subalgebra generated by these entries. Substituting $X$ and the $Y_{ij}$ for $1+r_{1}r_{2}$ independent matrices gives a surjective homomorphism from their polynomial ring onto $B$. This map is $G$-equivariant for simultaneous conjugation, since $(gv_{j})(u_{i}g^{-1})=gY_{ij}g^{-1}$. By Proposition \ref{p2.2}, it is still surjective on invariants. By the first fundamental theorem for matrices in Section 2.5 of \cite{c12}, the invariants of the independent matrices are generated by the traces of their products. Hence every invariant in $B$ is a polynomial in traces of products of $X$ and the $Y_{ij}$.
		\par A product involving only $X$ has trace $tr(X^{q})$ for some $q\ge 1$. For a product containing at least one $Y_{ij}$, cyclically permute its factors and substitute $Y_{ij}=v_{j}u_{i}$. Row and column multiplication gives
		\begin{align*}
			tr\big(Y_{i_{1}j_{1}}X^{l_{1}}\cdots Y_{i_{s}j_{s}}X^{l_{s}}\big)
			=\prod_{h=1}^{s}u_{i_{h}}(X^{l_{h}}v_{j_{h+1}}),
		\end{align*}
		where $j_{s+1}=j_{1}$ and $l_{h}\ge 0$. This proves that the stated functions generate all invariants. Conversely, each of these functions is invariant under the given action.
	\end{proof}

	\section{Proof of the main theorem}\label{s3}
	
	\par We prove the main theorem in this section.
	
	\par The determinantal variety $\mathcal{L}^{m,n}_{1}$ is one point for all $m,n \in \mathbb{Z}_{\ge 0}$. The main theorem is trivial in these cases. Throughout this section, suppose $r, m, n\in \mathbb{Z}_{\ge 2}$ are fixed integers satisfying $r\le m\le n$, and suppose $k\in \mathbb{Z}_{\ge 0}$.
	
	\par Every matrix $X\in \mathcal{L}_{r}^{m,n}$ admits a factorization $X=AB$, where $A\in M(m\times (r-1),K)$ and $B\in M((r-1)\times n,K)$. Suppose $\mathcal{U}^{m,n}_{r}=M(m\times (r-1),K)\times M((r-1)\times n,K)$. Then we have a surjective map
	\begin{align*}
		\mu:\mathcal{U}^{m,n}_{r}\to \mathcal{L}_{r}^{m,n},(A,B)\mapsto AB.
	\end{align*}
	At the level of $k$-jets, we write $X(t)=\sum\limits_{l=0}^{k}X^{(l)}t^{l}$ with $X^{(l)}=(x_{ij}^{(l)})_{1\le i\le m,1\le j\le n}$ for all $0\le l\le k$. Similarly, we have the notations $A(t)$ and $B(t)$ for $\mathcal{U}_{r}^{m,n}$. We write $K[J_{k}(\mathcal{U}_{r}^{m,n})]=K[a,b]$ for convenience. The map $J_{k}(\mu)$ is induced by $X(t)=A(t)B(t)$. Hence we have
	\begin{align*}
		X^{(l)}=\sum\limits_{l_{1}+l_{2}=l}A^{(l_{1})}B^{(l_{2})}
	\end{align*}
	for all $0\le l\le k$. Suppose $u_{ij}^{(l)}=\sum\limits_{l_{1}+l_{2}=l}\sum\limits_{h=1}^{r-1}a_{ih}^{(l_{1})}b_{hj}^{(l_{2})}\in K[a,b]$ for all $0\le l\le k$. By Corollary \ref{c2.1}, $K[P_{k}(\mathcal{L}_{r}^{m,n})]$ is isomorphic to 
	\begin{align*}
		K[u_{ij}^{(l)}|1\le i\le m,1\le j\le n,0\le l\le k],
	\end{align*}
	which is denoted by $D_{0}$ in the rest of the paper. By Proposition \ref{p2.1}, we have $H_{P_{k}(\mathcal{L}_{r}^{m,n})}(z)=H_{D_{0}}(z)$.
	
	\par $\quad$
	\par \noindent \textbf{Step 1. One-parameter deformation of $J_{k}(\mathcal{U}^{m,n}_{r})$.}
	
	\par Choose pairwise distinct $\lambda_{0},\cdots ,\lambda_{k}\in K$. We introduce the polynomial $f_{\varepsilon}(t)=(t-\lambda_{0}\varepsilon)\cdots (t-\lambda_{k}\varepsilon)\in K[\varepsilon,t]$. For all $1\le i\le m,1\le j\le n$, there exists a unique family of polynomials $\{ c_{ij}^{(l)}(\varepsilon)\in K[\varepsilon,a,b]|0\le l\le k\}$, such that the $(i,j)$-entry of $A(t)B(t)$ satisfies
	\begin{align*}
		\big( A(t)B(t)\big)_{ij}\equiv \sum\limits_{l=0}^{k}c_{ij}^{(l)}(\varepsilon)t^{l} \text{ (mod $f_{\varepsilon}(t)$)}.
	\end{align*}
	We define $C=K[\varepsilon][c_{ij}^{(l)}(\varepsilon)|1\le i\le m,1\le j\le n,0\le l\le k]$, which is a $K[\varepsilon]$-subalgebra of $K[\varepsilon,a,b]$. Evaluation at $\varepsilon=0$ induces $ev_{\varepsilon=0}:K[\varepsilon,a,b]\to K[a,b]$.
	
	\begin{lemma}\label{l3.1}
		For all $1\le i\le m,1\le j\le n, 0\le l\le k$, we have $ev_{\varepsilon=0}(c_{ij}^{(l)}(\varepsilon))=u_{ij}^{(l)}$.
	\end{lemma}
	\begin{proof}[Proof. ]
		Choose a polynomial $\varphi(t)\in K[\varepsilon,a,b][t]$, such that
		\begin{align*}
			\big( A(t)B(t)\big)_{ij}=\sum\limits_{l=0}^{k}c_{ij}^{(l)}(\varepsilon)t^{l}+f_{\varepsilon}(t)\varphi(t).
		\end{align*}
		Taking $\varepsilon=0$, we obtain
		\begin{align*}
			\big( A(t)B(t)\big)_{ij}=\sum\limits_{l=0}^{k}ev_{\varepsilon=0}(c_{ij}^{(l)}(\varepsilon))t^{l}+t^{k+1}\varphi(t),
		\end{align*}	
		which completes the proof.
	\end{proof}
	
	\par We can restrict $ev_{\varepsilon=0}$ to $C$. The lemma above shows that $ev_{\varepsilon=0}(C)\subset D_{0}$. By definition, $ev_{\varepsilon=0}(\varepsilon C)=0$, so we can define $ev_{\varepsilon=0}:C/\varepsilon C\to D_{0}$.
	
	\begin{lemma}\label{l3.2}
		The $K$-algebra homomorphism $ev_{\varepsilon=0}:C/\varepsilon C\to D_{0}$ is surjective.
	\end{lemma}
	\begin{proof}[Proof. ]
		This is a direct corollary of Lemma \ref{l3.1}.
	\end{proof}
	
	\begin{remark}\label{r3.1}
		In fact, $ev_{\varepsilon=0}:C/\varepsilon C\to D_{0}$ is an isomorphism. This is not easy to prove and we will prove it later. 
	\end{remark}

	\par $\quad$
	\par \noindent \textbf{Step 2. $H_{C/\varepsilon C}(z)=H_{\mathcal{L}_{r}^{m,n}}(z)^{k+1}$.}
	
	\par The $K[\varepsilon]$-module $K[\varepsilon,a,b]$ is a free $K[\varepsilon]$-module. The ring $K[\varepsilon]$ is a PID. Hence we have $H_{C/\varepsilon C}(z)=H_{C\otimes_{K[\varepsilon]}K(\varepsilon),K(\varepsilon)}(z)$. For any two finitely generated graded $K$-algebras $R_{1}$ and $R_{2}$, we have $H_{R_{1}\otimes_{K}R_{2}}(z)=H_{R_{1}}(z)H_{R_{2}}(z)$. Hence $H_{\mathcal{L}_{r}^{m,n}}(z)^{k+1}=H_{K[\mathcal{L}_{r}^{m,n}]^{\otimes_{K}(k+1)}}(z)$. The Hilbert series is invariant under the field extension, so we have
	\begin{align*}
		H_{\mathcal{L}_{r}^{m,n}}(z)^{k+1}=H_{K[\mathcal{L}_{r}^{m,n}]^{\otimes_{K}(k+1)}\otimes_{K}K(\varepsilon),K(\varepsilon)}(z).%=H_{(K[\mathcal{L}_{r}^{m,n}]\otimes_{K}K(\varepsilon))^{\otimes_{K(\varepsilon)}(k+1)},K(\varepsilon)}(z).
	\end{align*}
	We write $D^{\prime}=K[\mathcal{L}_{r}^{m,n}]^{\otimes_{K}(k+1)}$ for convenience.
	
	\par For all $0\le \alpha\le k$, suppose $\varphi_{\alpha}:K[\mathcal{L}_{r}^{m,n}]\to D^{\prime}$ is the $(\alpha+1)$th factor embedding. Consider the $K(\varepsilon)$-linear map
	\begin{align*}
		\tau:D^{\prime}\otimes_{K}K(\varepsilon)\to C\otimes_{K[\varepsilon]}K(\varepsilon), \varphi_{\alpha}(x_{ij})\otimes 1\mapsto \sum\limits_{l=0}^{k}(\lambda_{\alpha}\varepsilon)^{l}c_{ij}^{(l)}(\varepsilon).
	\end{align*}
	We have to show that $\tau$ is an isomorphism for this step.
	
	\begin{prop}\label{p3.1}
	    Notations are as above, $\tau:D^{\prime}\otimes_{K}K(\varepsilon)\to C\otimes_{K[\varepsilon]}K(\varepsilon)$ is an isomorphism.
	\end{prop}
	\begin{proof}[Proof. ]
		For all $0\le \alpha\le k$, we take $t=\lambda_{\alpha}\varepsilon$ for $A(t)$ and $B(t)$, and suppose
		\begin{center}
			$a_{ij,\alpha}=A(\lambda_{\alpha}\varepsilon)_{ij}=\sum\limits_{l=0}^{k}a_{ij}^{(l)}(\lambda_{\alpha}\varepsilon)^{l}$,\\
			$b_{ij,\alpha}=B(\lambda_{\alpha}\varepsilon)_{ij}=\sum\limits_{l=0}^{k}b_{ij}^{(l)}(\lambda_{\alpha}\varepsilon)^{l}$
		\end{center}
		for all $i,j$. Let $V_{\varepsilon}=\big((\lambda_{\alpha}\varepsilon)^{l}\big)_{0\le \alpha\le k,0\le l\le k}$ be a $(k+1)\times (k+1)$ matrix. We have
		\begin{align*}
			\begin{pmatrix}
				a_{ij,0}\\
				\vdots \\
				a_{ij,k}
			\end{pmatrix}
			=V_{\varepsilon}
			\begin{pmatrix}
			    a_{ij}^{(0)}\\
			    \vdots \\
			    a_{ij}^{(k)}
			\end{pmatrix},\begin{pmatrix}
			b_{ij,0}\\
			\vdots \\
			b_{ij,k}
			\end{pmatrix}
			=V_{\varepsilon}
			\begin{pmatrix}
			b_{ij}^{(0)}\\
			\vdots \\
			b_{ij}^{(k)}
			\end{pmatrix},
		\end{align*}
		and the determinant of $V_{\varepsilon}$ is 
		\begin{align*}
			\varepsilon^{\frac{k(k+1)}{2}}\prod\limits_{0\le \beta<\alpha\le k}(\lambda_{\alpha}-\lambda_{\beta})\neq 0.
		\end{align*}
		Hence we obtain an isomorphism
		\begin{align*}
			K[a,b]\simeq K[a_{ih,\alpha},b_{hj,\alpha}|1\le i\le m,1\le j\le n,1\le h\le r-1,0\le \alpha\le k].
		\end{align*}
		For all $0\le \alpha\le k$, let $\widehat{C}_{\alpha}=K[a_{ih,\alpha},b_{hj,\alpha}|1\le i\le m,1\le j\le n,1\le h\le r-1]$. The isomorphism above is $K[a,b]\simeq \widehat{C}_{0}\otimes \cdots \otimes \widehat{C}_{k}$, and $\tau \circ \varphi_{\alpha}:K[\mathcal{L}_{r}^{m,n}]\otimes_{K}K(\varepsilon)\to \widehat{C}_{\alpha}$ is injective. It follows that $\tau$ is injective.
		\par For all $0\le \alpha\le k$, we write $\widehat{A}_{\alpha}=(a_{ih,\alpha})_{1\le i\le m,1\le h\le r-1}$ and $\widehat{B}_{\alpha}=(b_{hj,\alpha})_{1\le h\le r-1,1\le j\le n}$.\\\\
		\textit{Claim}. For all $i,j,\alpha$, we have $\big(\widehat{A}_{\alpha}\widehat{B}_{\alpha}\big)_{ij}=\sum\limits_{l=0}^{k}(\lambda_{\alpha}\varepsilon)^{l}c_{ij}^{(l)}(\varepsilon)$.
		\begin{proof}[Proof of Claim. ]
			There exists a polynomial $\varphi_{ij}(t)\in K[\varepsilon,a,b][t]$ such that
			\begin{align*}
				\big(A(t)B(t)\big)_{ij}=\sum\limits_{l=0}^{k}c_{ij}^{(l)}(\varepsilon)t^{l}
				+f_{\varepsilon}(t)\varphi_{ij}(t).
			\end{align*}
			Taking $t=\lambda_{\alpha}\varepsilon$, we have $\big(A(\lambda_{\alpha}\varepsilon)B(\lambda_{\alpha}\varepsilon)\big)_{ij}-\sum\limits_{l=0}^{k}c_{ij}^{(l)}(\varepsilon)(\lambda_{\alpha}\varepsilon)^{l}=f_{\varepsilon}(\lambda_{\alpha}\varepsilon)\varphi_{ij}(\lambda_{\alpha}\varepsilon)=0$. Since $A(\lambda_{\alpha}\varepsilon)=\widehat{A}_{\alpha}$ and $B(\lambda_{\alpha}\varepsilon)=\widehat{B}_{\alpha}$, this gives the claimed equality.
		\end{proof}
		\par By the claim, we have
		\begin{align*}
			\tau\circ \varphi_{\alpha}:K[\mathcal{L}_{r}^{m,n}]\otimes_{K}K(\varepsilon)\longrightarrow\widehat{C}_{\alpha},
			x_{ij}\otimes 1\longmapsto\big(\widehat{A}_{\alpha}\widehat{B}_{\alpha}\big)_{ij}.
		\end{align*}
		We define $L_{\alpha}(t)=\prod\limits_{\substack{0\le\beta\le k,\beta\neq \alpha}}\frac{t-\lambda_{\beta}\varepsilon}{\lambda_{\alpha}\varepsilon-\lambda_{\beta}\varepsilon}\in K(\varepsilon)[t]$. The $\lambda_{\alpha}$ are distinct, so $L_{\alpha}(\lambda_{\beta}\varepsilon)=\delta_{\alpha\beta}$. Hence we have
		\begin{align*}
			\sum\limits_{l=0}^{k}c_{ij}^{(l)}(\varepsilon)t^{l}
			=\sum\limits_{\alpha=0}^{k}\tau\circ \varphi_{\alpha}(x_{ij}\otimes 1)L_{\alpha}(t),
		\end{align*}
		since both sides have degree at most $k$ in $t$ and agree at $\lambda_{0}\varepsilon,\cdots,\lambda_{k}\varepsilon$. Comparing coefficients of $t^{l}$, one can see that every $c_{ij}^{(l)}(\varepsilon)$ is in $\Ima \tau$. So $\tau$ is surjective.
	\end{proof}
	
	\par By Proposition \ref{p2.1} and Lemma \ref{l3.2}, it remains to show $ev_{\varepsilon=0}:C/\varepsilon C\to D_{0}$ is injective.

	\par $\quad$
	\par \noindent \textbf{Step 3. Flag variety structure of jet schemes and conjugate action.}
	
	\par Let $V=\big( K[t]/(t^{k+1})\big)^{r-1}$, and let $\Gamma=K[\varepsilon,t]/(f_{\varepsilon}(t))$. We have $\Gamma^{r-1}\simeq K[\varepsilon]\otimes_{K}V$. Consider $\{ t^{l}e_{h}|0\le l\le k,1\le h\le r-1\}$, where $e_{h}=(0,\cdots ,0,1,0,\cdots ,0)$. It is a $K$-basis for $V$ and a $K[\varepsilon]$-basis for $\Gamma^{r-1}$. We will use the flag variety $F=F_{r-1,\cdots ,r-1}(V)$, where $r-1$ occurs $k+1$ times. Notations about flag variety in Section 2.3 are kept. 
	\par Multiplication by $t$ defines a $K[\varepsilon]$-linear map
	\begin{align*}
		T:\Gamma^{r-1}\to \Gamma^{r-1},x\mapsto tx.
	\end{align*}
	For all $c\in K$, we denote by $T_{c}\in \End(V)$ the specialization of $T$ at $\varepsilon=c$, using the basis above. Thus $T$ is represented by a matrix with entries in $K[\varepsilon]$, while every $T_{c}$ is a linear map on the finite dimensional $K$-vector space $V$. In particular, $T_{0}$ is multiplication by $t$ on $V$. One can see that the Jordan form of $T_{0}$ consists of $r-1$ nilpotent Jordan blocks of size $k+1$. Hence the nilpotent orbit is $\mathcal{N}=\overline{\{ gT_{0}g^{-1}|g\in G\}}\subset \End(V)$.
	
	\begin{lemma}\label{l3.4}
		For all $c\in K$ nonzero, $T_{c}$ is conjugate to
		\begin{align*}
			\operatorname{diag}(\lambda_{0}cI_{r-1},\cdots ,\lambda_{k}cI_{r-1}).
		\end{align*}
	\end{lemma}
	\begin{proof}[Proof. ]
		Write $f_{c}(t)=f_{\varepsilon}(t)|_{\varepsilon=c}$. Since $\lambda_{0}c,\cdots ,\lambda_{k}c$ are pairwise distinct, evaluation gives an isomorphism
		\begin{align*}
			K[t]/(f_{c}(t))\simeq K^{k+1},p(t)\mapsto \big( p(\lambda_{0}c),\cdots ,p(\lambda_{k}c)\big).
		\end{align*}
		Multiplication by $t$ acts on the $(\alpha+1)$th factor as multiplication by $\lambda_{\alpha}c$. Applying this isomorphism to all $r-1$ factors gives the result.
	\end{proof}
	
	\par The variety $Z$ in Section 2.3 for this $F$ is 
	\begin{align*}
		Z=\{ (D,F_{\bullet},\varepsilon)\in \End(V)\times F\times \mathbb{A}_{K}^{1}|(D-\lambda_{i}\varepsilon I)F_{i}\subset F_{i+1} \text{ for all $0\le i\le k$}\}.
	\end{align*}
	By Section 2.3, $Z_{0}\simeq T^{*}F$ is the nilpotent orbit $\mathcal{N}$. The map $\mu:Z\to \End(V)\times \mathbb{A}_{K}^{1}$ is projective by Proposition \ref{p2.4}. The group $G$ acts on $Z$ by
	\begin{align*}
		g\cdot(D,F_{\bullet},\varepsilon)=(gDg^{-1},gF_{\bullet},\varepsilon).
	\end{align*}
	The projection $\mu$ is equivariant for this action and the action $g\cdot(D,\varepsilon)=(gDg^{-1},\varepsilon)$ on its target.
	
	\par We now put the nonzero fibers into one affine variety. Define
	\begin{align*}
		Y=\overline{\{ (gT_{c}g^{-1},c)|g\in G,c\in K^{*}\}}\subset \End(V)\times \mathbb{A}_{K}^{1}.
	\end{align*}
	The projection $\mu$ induces a projective surjective map $\mu:Z\to Y$.
	
	\begin{lemma}\label{l3.6}
		The map $\mu:Z\to Y$ restricts to an isomorphism over $\mathbb{A}_{K}^{1}\setminus\{0\}$.
	\end{lemma}
	\begin{proof}[Proof. ]
		For all $c\in K$ nonzero and $D\in\{ gT_{c}g^{-1}|g\in G\}$, we have $\mu(F_{\bullet})=D$, where
		\begin{align*}
			F_{i}=\bigoplus\limits_{\alpha=i}^{k}\ker(D-\lambda_{\alpha}cI)
		\end{align*}
		for all $0\le i\le k+1$. Hence $\mu$ is bijective. It remains to show that the inverse constructed as above varies regularly with $c$. Let
		\begin{align*}
			P_{\alpha}(D,c)=\prod\limits_{\substack{0\le\beta\le k,\beta\neq\alpha}}\frac{D-\lambda_{\beta}cI}{(\lambda_{\alpha}-\lambda_{\beta})c}.
		\end{align*}
		They are regular for $c\neq 0$, and $\Ima P_{\alpha}(D,c)=\ker(D-\lambda_{\alpha}cI)$ has constant dimension $r-1$. Thus the subspaces $F_{i}$ define a regular map to $F$.
	\end{proof}

	\par $\quad$
	\par \noindent \textbf{Step 4. Group action on $\mathcal{U}_{r}^{m,n}$ via base change.}
	
	\par Suppose $H=GL(r-1,K)$, and we consider the $H$-action on $\mathcal{U}_{r}^{m,n}$
	\begin{align*}
		q\cdot (A,B)=(Aq^{-1},qB).
	\end{align*}
	We will use this group action and Proposition \ref{p2.3} to construct relations between  $\mathcal{U}_{r}^{m,n}$ and representation of general linear groups.
	
	\par We identify a row vector $(a_{1},\cdots ,a_{r-1})$ with $r-1$ entries in $\Gamma$ with an element of $K[\varepsilon]\otimes_{K}V^{*}$, via the standard dot-product pairing. Hence we obtain an isomorphism
	\begin{align*}
		\Psi:M(m\times (r-1),\Gamma)\times M((r-1)\times n,\Gamma)\to K[\varepsilon]\otimes_{K}W,
	\end{align*}
	where $W=(V^{*})^{m}\oplus V^{n}$. For $f\in K[\varepsilon,t]$, suppose $f$ is written as $\widetilde{f}$ after modulo $(f_{\varepsilon}(t))$, where $\widetilde{f}$ is a polynomial of degree less that $k+1$ in $t$. Let $\eta(f)=[t^{k}](\widetilde{f})$, where $[t^{k}]$ denotes the operation of taking the coefficient of $t^{k}$. More precisely, $\Psi(A(t),B(t))=(u_{1},\cdots ,u_{m},v_{1},\cdots ,v_{n})$, where
	\begin{equation}\label{eq-uv}
		u_{i}(t^{l}e_{h})=\eta(a_{ih}(t)t^{l}), v_{j}=b_{1j}(t)e_{1}+\cdots +b_{r-1,j}(t)e_{r-1}.
	\end{equation}
	Its specialization at $\varepsilon=0$ is denoted by $\Psi_{0}:J_{k}(\mathcal{U}_{r}^{m,n})\to W$. For $\Psi_{0}$, we have $u_{i}(t^{l}e_{h})=a_{ih}^{(k-l)}$ for all $0\le l\le k$.
	
	\par The same formula can be used over $\Gamma$. For all $q\in GL(r-1,\Gamma)$, we have
	\begin{align*}
		\Psi(Aq^{-1},qB)=(u_{1}\circ q^{-1},\cdots ,u_{m}\circ q^{-1},qv_{1},\cdots ,qv_{n}).
	\end{align*}
	At $\varepsilon=0$, the formula applies to every $q\in GL(r-1,K[t]/(t^{k+1}))$ and agrees with the $J_{k}(H)$-action on $J_{k}(\mathcal{U}_{r}^{m,n})$. Since these changes of basis commute with multiplication by $t$, the formula agrees with the action in Proposition \ref{p2.3} at a fixed operator.
	
	\par We keep the $G$-action on $Y$ and $\mathcal{N}$ in Step 3. Together with the action on $W$ in Proposition \ref{p2.3}, it gives the diagonal $G$-action on $Y\times W$ and $\mathcal{N}\times W$, with $\varepsilon$ fixed. We write
	\begin{align*}
		\mathcal{A}=K[Y\times W]^{G},\mathcal{A}_{0}=K[\mathcal{N}\times W]^{G}.
	\end{align*}
	For all $1\le i\le m$, $1\le j\le n$, and $l\ge0$, let $w_{ij}^{(l)}=u_{i}(D^{l}v_{j})$ be the corresponding invariant function on either variety.
	
	\begin{lemma}\label{l3.8}
		We have $\mathcal{A}=K[\varepsilon][w_{ij}^{(l)}|1\le i\le m,1\le j\le n,0\le l\le k]$.
	\end{lemma}
	\begin{proof}[Proof. ]
		The restriction from $K[\End(V)\oplus W][\varepsilon]$ to $K[Y\times W]$ is surjective and $G$-equivariant. By Propositions \ref{p2.2} and \ref{p2.3}, we obtain
		\begin{align*}
			\mathcal{A}=K[\varepsilon][tr(D^{q}),w_{ij}^{(l)}|q\ge1,l\ge0,1\le i\le m,1\le j\le n].
		\end{align*}
		By Lemma \ref{l3.4}, on the dense family defining $Y$, and hence on $Y$, we have
		\begin{align*}
			f_{\varepsilon}(D)=0,tr(D^{q})=(r-1)\sum\limits_{\alpha=0}^{k}(\lambda_{\alpha}\varepsilon)^{q}\text{ for all $q\ge1$}.
		\end{align*}
		Thus the traces belong to $K[\varepsilon]$. We have the relation $f_{\varepsilon}(D)=0$. Using this relation, we can write every power $D^{l}$ with $l\ge k+1$ as a $K[\varepsilon]$-linear combination of $D^{0},\cdots ,D^{k}$.
	\end{proof}
	
	\begin{lemma}\label{l3.9}
		We have $\mathcal{A}_{0}=K[w_{ij}^{(l)}|1\le i\le m,1\le j\le n,0\le l\le k]$.
	\end{lemma}
	\begin{proof}[Proof. ]
		Apply Propositions \ref{p2.2} and \ref{p2.3} to the restriction from $K[\End(V)\oplus W]$ to $K[\mathcal{N}\times W]$. The invariant ring is generated by $tr(D^{q})$ and $w_{ij}^{(l)}$. Since $T_{0}^{k+1}=0$, we have $D^{k+1}=0$ and $tr(D^{q})=0$ on $\mathcal{N}$ for all $q\ge 1$. The result follows.
	\end{proof}
	
	\par The section $c\mapsto(T_{c},c)$ takes values in $Y$ for all $c\in K$. Indeed, this holds for $c\neq0$ by the definition of $Y$, and then for $c=0$ since $Y$ is closed. Restriction to this section, followed by the pullback under $\Psi$, defines a $K[\varepsilon]$-algebra homomorphism
	\begin{align*}
		\theta:\mathcal{A}\to K[\varepsilon,a,b],f\mapsto\Psi^{*}\big(f(T,\varepsilon,u,v)\big).
	\end{align*}
	Here we follow notations (\ref{eq-uv}) with $u=(u_{1},\cdots ,u_{m})$ and $v=(v_{1},\cdots ,v_{n})$. The expression $f(T,\varepsilon,u,v)$ is the polynomial function on $\mathbb{A}_{K}^{1}\times W$ whose value at $\varepsilon=c$ is $f(T_{c},c,u,v)$. Similarly, restriction to $T_{0}\in\mathcal{N}$, followed by the pullback under $\Psi_{0}$, defines a $K$-algebra homomorphism
	\begin{align*}
		\theta_{0}:\mathcal{A}_{0}\to K[a,b],f\mapsto\Psi_{0}^{*}\big(f(T_{0},u,v)\big).
	\end{align*}
	
	\par We now calculate their images on the generators. We $K[\varepsilon,a,b]$-linearly extend $\eta$ when the coefficient variables $a,b$ are present. By the definitions of $\Psi$ and $c_{ij}^{(a)}(\varepsilon)$, we have
	\begin{align*}
		\theta(w_{ij}^{(l)})=\eta\big(t^{l}(A(t)B(t))_{ij}\big)=\sum\limits_{a=0}^{k}\eta(t^{l+a})c_{ij}^{(a)}(\varepsilon)\in C
	\end{align*}
	for all $0\le l\le k$. Since $\theta$ commutes with $\varepsilon$, Lemma \ref{l3.8} shows that $\theta(\mathcal{A})\subset C$. At $\varepsilon=0$, reduction is modulo $t^{k+1}$. Therefore, taking the coefficient of $t^{k}$ gives $\theta_{0}(w_{ij}^{(l)})=u_{ij}^{(k-l)}\text{ for all $0\le l\le k$}$. By Lemma \ref{l3.9}, $\theta_{0}(\mathcal{A}_{0})\subset D_{0}$. Hence we obtain the two homomorphisms
	\begin{align*}
		\theta:\mathcal{A}\to C,\theta_{0}:\mathcal{A}_{0}\to D_{0}.
	\end{align*}
	We will prove that these maps are isomorphisms later.

	\par $\quad$
	\par \noindent \textbf{Step 5. $\mathcal{A}/\varepsilon \mathcal{A}$ and $\mathcal{A}_{0}$ are isomorphic.}
	
	\par We will abuse notation by using the same symbol for a coherent sheaf and its global sections.
	
	\par We have $K[\mathcal{N}]\simeq K[Z_{0}]$, as discussed in the previous steps. The closed embedding $i:Z_{0}\to Z$ gives a short exact sequence
	\begin{equation}\label{eq5}
		0\to \varepsilon \mathcal{O}_{Z}\to \mathcal{O}_{Z}\to i_{*}\mathcal{O}_{Z_{0}}\to 0.
	\end{equation}
	
	\begin{lemma}\label{l3.3}
		The map $\mu:Z\to Y$ induces an isomorphism $K[Y]\simeq K[Z]$.
	\end{lemma}
	\begin{proof}[Proof. ]
		Since $\mu$ is surjective, the pullback $\mu^{*}:K[Y]\to K[Z]$ is injective. We identify $K[Y]$ with its image in $K[Z]$. Since $\mu$ is projective and $Y$ is affine, $K[Z]$ is a finitely generated $K[Y]$-module. Hence $L=K[Z]/K[Y]$ is also a finitely generated $K[Y]$-module. Restriction to $Z_{0}$ induces a surjective $i^{*}\mu^{*}:K[Y]\to K[Z_{0}]\simeq K[Z]/\varepsilon K[Z]$. For any $f\in K[Z]$, the surjectivity of $i^{*}\mu^{*}$ gives $g\in K[Y]$ such that $i^{*}(f)=i^{*}\mu^{*}(g)$. Hence $f-\mu^{*}(g)\in\ker i^{*}=\varepsilon K[Z]$. Thus we have $K[Z]=K[Y]+\varepsilon K[Z]$, and therefore $L=\varepsilon L$. By Lemma \ref{l3.6}, we have
		\begin{align*}
			K[Y][\varepsilon^{-1}]=K[Z][\varepsilon^{-1}].
		\end{align*}
		Hence $L[\varepsilon^{-1}]=0$. Since $L$ is finitely generated over $K[Y]$, there exists an integer $s\ge 0$ such that $\varepsilon^{s}L=0$. Combining this with $L=\varepsilon L$, we complete the proof.
	\end{proof}
	
	\par The short exact sequence (\ref{eq5}) gives a short exact sequence
	\begin{align*}
		0\to \varepsilon K[Y]\to K[Y]\to K[\mathcal{N}]\to 0
	\end{align*}
	via the lemma above. By the Künneth formula, we have a short exact sequence
	\begin{align*}
		0\to \varepsilon K[Y\times W]\to K[Y\times W]\to K[\mathcal{N}\times W]\to 0
	\end{align*}
	since $W$ is an affine space. Applying Proposition \ref{p2.2} to this exact sequence, we obtain an isomorphism
	\begin{align*}
		\sigma:K[Y\times W]^{G}/(\varepsilon K[Y\times W])^{G}\to K[\mathcal{N}\times W]^{G}.
	\end{align*}
	The $G$-action commutes with $\varepsilon$, so we can write $\sigma:\mathcal{A}/\varepsilon \mathcal{A}\to \mathcal{A}_{0}$. Since $\theta$ commutes with $\varepsilon$, $\theta$ induces a homomorphism $\widetilde{\theta}:\mathcal{A}/\varepsilon \mathcal{A}\to C/\varepsilon C$.
	
	\begin{lemma}\label{lcom}
		The diagram
		\begin{figure}[H]
			\centering
			\begin{tikzpicture}
				\node(a) at (0,2){$\mathcal{A}/\varepsilon \mathcal{A}$};
				\node(b) at (3,2){$\mathcal{A}_{0}$};
				\node(c) at (0,0){$C/\varepsilon C$};   
				\node(d) at (3,0){$D_{0}$};
				\draw[->] (a)--(b) node[midway,above]{$\sigma$};
				\draw[->] (b)--(d) node[midway,right]{$\theta_{0}$};
				\draw[->] (a)--(c) node[midway,left]{$\widetilde{\theta}$};
				\draw[->] (c)--(d) node[midway,above]{$ev_{\varepsilon =0}$};
			\end{tikzpicture}
		\end{figure}	\FloatBarrier
		\par \noindent is commutative.
	\end{lemma}
	\begin{proof}[Proof. ]
		By Lemma \ref{l3.8}, $\mathcal{A}/\varepsilon\mathcal{A}$ is generated by the classes $\overline{w_{ij}^{(l)}}$ for all $1\le i\le m$, $1\le j\le n$, and $0\le l\le k$. Since $\sigma(\overline{w_{ij}^{(l)}})=w_{ij}^{(l)}$, we have $\theta_{0}\sigma(\overline{w_{ij}^{(l)}})=u_{ij}^{(k-l)}$. On the other hand, at $\varepsilon=0$, reduction is modulo $t^{k+1}$, so
		\begin{align*}
			\eta(t^{l+a})\big|_{\varepsilon=0}
			=\begin{cases}
				1,&l+a=k,\\
				0,&l+a\neq k.
			\end{cases}
		\end{align*}
		By Lemma \ref{l3.1}, it follows that
		\begin{align*}
			ev_{\varepsilon=0}\widetilde{\theta}(\overline{w_{ij}^{(l)}})=\sum\limits_{a=0}^{k}\eta(t^{l+a})\big|_{\varepsilon=0}u_{ij}^{(a)}=u_{ij}^{(k-l)}.
		\end{align*}
		The two compositions agree on all generators, hence the diagram is commutative.
	\end{proof}

	\par $\quad$
	\par \noindent \textbf{Step 6. $\theta$ and $\theta_{0}$ are isomorphisms.}
	
	\par Since $\theta$ commutes with $\varepsilon$, it is a direct corollary of this step that $\widetilde{\theta}$ is an isomorphism. Then it follows from Lemma \ref{lcom} that $ev_{\varepsilon=0}$ is an isomorphism, which completes the proof of the main theorem.
	
	\begin{lemma}\label{1isom}
		The homomorphism $\theta_{0}:\mathcal{A}_{0}\to D_{0}$ is an isomorphism.
	\end{lemma}
	\begin{proof}[Proof. ]
		We have $\theta_{0}(w_{ij}^{(l)})=u_{ij}^{(k-l)}$ for all $0\le l\le k$. The image of $\theta_{0}$ contains every generator of $D_{0}$, so $\theta_{0}$ is surjective.
		
		\par Suppose $f\in\mathcal{A}_{0}$ and $\theta_{0}(f)=0$. At $\varepsilon=0$, the coordinates in (\ref{eq-uv}) satisfy $u_i(t^le_h)=a_{ih}^{(k-l)}$, and the coefficients of the $v_j$ are the coefficients of $B(t)$. Thus $\Psi_{0}$ is an invertible change of coordinates. By the definition of $\theta_{0}$, we obtain $f(T_{0},w)=0$ for all $w\in W$. For every $g\in G$, the $G$-invariance of $f$ gives
		\begin{align*}
			f(gT_{0}g^{-1},w)=f(T_{0},g^{-1}\cdot w)=0.
		\end{align*}
		The subset $(G\cdot T_{0})\times W$ is dense in the variety defining $\mathcal{A}_{0}$. Therefore we have $f=0$, and $\theta_{0}$ is injective.
	\end{proof}
	
	\begin{lemma}\label{2isom}
		The homomorphism $\theta:\mathcal{A}\to C$ is an isomorphism.
	\end{lemma}
	\begin{proof}[Proof. ]
		For all $0\le l,p\le k$, the definition of $\eta$ gives
		\begin{align*}
			\eta(t^{l+p})=0\text{ if $l+p<k$},\\
			\eta(t^{l+p})=1\text{ if $l+p=k$}.
		\end{align*}
		After reversing the columns, the matrix $\big(\eta(t^{l+p})\big)_{0\le l,p\le k}$ is lower triangular with diagonal entries one. Hence this matrix is invertible over $K[\varepsilon]$. For fixed $i,h$, (\ref{eq-uv}) gives
		\begin{align*}
			u_i(t^le_h)=\sum\limits_{p=0}^{k}\eta(t^{l+p})a_{ih}^{(p)}.
		\end{align*}
		Thus the row part of $\Psi$ is an invertible change of coefficients over $K[\varepsilon]$. The coefficients of $B(t)$ are exactly the coordinates of the $v_j$ in the basis $t^le_h$, so $\Psi$ is invertible over $K[\varepsilon]$.
		
		\par Suppose $f\in\mathcal{A}$ and $\theta(f)=0$. Since $\Psi$ is invertible, the definition of $\theta$ implies $f(T_c,c,w)=0$ for all $c\in K$ and $w\in W$. For every $g\in G$, the $G$-invariance of $f$ gives
		\begin{align*}
			f(gT_cg^{-1},c,w)=f(T_c,c,g^{-1}\cdot w)=0.
		\end{align*}
		By the definition of $Y$, these points with $c\neq 0$ form a dense subset of $Y\times W$. Therefore we have $f=0$, and $\theta$ is injective.
		
		\par To prove surjectivity, recall that $\theta(\varepsilon)=\varepsilon$ and
		\begin{align*}
			\theta(w_{ij}^{(l)})=\sum\limits_{a=0}^{k}\eta(t^{l+a})c_{ij}^{(a)}(\varepsilon)
		\end{align*}
		for all $0\le l\le k$. For fixed $i,j$, the coefficient matrix of these equations is the invertible matrix considered above. Its inverse has entries in $K[\varepsilon]$, so every $c_{ij}^{(a)}(\varepsilon)$ is a $K[\varepsilon]$-linear combination of the $\theta(w_{ij}^{(l)})$. Therefore the image of $\theta$ contains all generators of $C$ over $K[\varepsilon]$, and $\theta$ is surjective.
	\end{proof}

	\section{Acknowledgement}
	
	AI tools (ChatGPT and DeepSeek) were used during the preparation of the paper as technical assistants to help with checking calculations and arguments, developing details of proofs, and improving the exposition. All mathematical statements, proofs, and conclusions were verified by the authors. We thank Qiyue Tang for helpful discussions with us. H. Zuo acknowledges support from NSFC (grant No. 12671056) and BJNSF (grant No. 1252009).
	
	\end{document}